\documentclass[a4paper,12pt]{amsart}
\usepackage[ansinew]{inputenc}
\usepackage[T1]{fontenc}

\usepackage{url}

\usepackage{amsmath}
\usepackage{amssymb}
\usepackage{amsfonts}
\usepackage{amsthm}

\usepackage{mathrsfs}

\newcommand{\Z}{\mathbb{Z}}
\newcommand{\Q}{\mathbb{Q}}

\newcommand{\F}{\mathbb{F}}

\newcommand{\imu}{\mathrm{i}}      
\newcommand{\neper}{\mathrm{e}}    

\newcommand{\Gal}{\text{Gal}}      
\newcommand{\Euler}{\varphi}       

\DeclareMathOperator{\Tr}{Tr}      
\DeclareMathOperator{\Kl}{Kl} 

\newcommand{\comment}[1]{}

\theoremstyle{plain}
  \newtheorem{theorem}{Theorem}
  \newtheorem{lemma}[theorem]{Lemma}
  \newtheorem{corollary}[theorem]{Corollary}
  
\title{On the degree of a Kloosterman sum as an algebraic integer}
  \author{Keijo Kononen}
  \email{keijo.kononen@oulu.fi}

  \author{Marko Rinta-aho}
  \email{marko.rinta-aho@oulu.fi}
  
  \author{Keijo Väänänen}
  \email{keijo.vaananen@oulu.fi}
  
  \address{Department of Mathematical Sciences, University Of Oulu, P.O. BOX 3000, FIN-90014 Oulun yliopisto, Finland}

\begin{document}

\begin{abstract}
The maximal degree over rational numbers that an $n$-dimensinonal Kloosterman sum defined over a finite field of characteristic $p$ can achieve is known to be $\tfrac{p-1}{d}$ where $d=\gcd(p-1,n+1)$. Wan has shown that this maximal degree is always achieved in points whose absolute trace is nonzero. By the works of Fischer, Wan we know that there exist many finite fields for which the values of the Kloosterman sums are distinct except Frobenius conjugation. For these fields we completely determine the degrees of all the Kloosterman sums. Even if the finite field does not satisfy this condition we can still often find points in which the Kloosterman sum has smaller degree than $\tfrac{p-1}{d}$. 
\end{abstract}
  \keywords{cyclotomic field, Kloosterman sum}

\maketitle

\section{Introduction}

Let $p,q,r$ be positive integers, $p$ a prime and $q=p^r$. We denote by $\F_{q}$ the finite field of $q$ elements and by $\Tr$ the absolute trace function 
$\F_{q}\to \F_p$. Also, let $\zeta_p:= \neper^{2\pi\imu/p}$ be a complex primitive $p$th root of unity. The $n$-dimensional Kloosterman sum at a point $a\in \F_{q}$ is defined by the equation 
\begin{equation*}
  \Kl_n(q,a) = \sum_{x_1,x_2,\ldots,x_n \in \F_q^*} \psi\left(x_1+x_2+\cdots+x_n +\frac{a}{x_1 x_2 \cdots x_n}\right), 
\end{equation*}  
where $\psi: \F_{q}\to \Q(\zeta_p)$ is the canonical additive character of the field $\F_{q}$ defined by $\psi(x)=\zeta_p^{\Tr(x)}$. 

Obviously the values of the Kloosterman sums are algebraic integers of the cyclotomic field $\Q(\zeta_p)$. For every index  $i\in\F_p^*$ we let $\sigma_i$ denote the unique $\Q$-automorphism of the field $\Q(\zeta_p)$ mapping $\zeta_p$ to $\zeta_p^i$. The mapping $i\mapsto \sigma_i$ is an isomorphism between the multiplicative group $\F_p^*$ of the prime field $\F_p$ and the Galois group $\Gal(\Q(\zeta_p)/\Q)$ of the extension $\Q(\zeta_p)/\Q$. Since $\F_p^*$ is cyclic there exists a unique subgroup of order $d$ for every positive factor $d$ of $p-1$. It follows by the main theorem of the Galois theory that there exists a unique intermediate field (denoted by $E_d$ from now on) of degree $\frac{p-1}{d}$ over $\Q$ between $\Q$ and $\Q(\zeta_p)$. This is the fixed field of the order $d$ subgroup. Since this is a cyclic group, $x\in E_d$ if and only if $\sigma_i(x)=x$ for some (and hence for all) $i\in\F_p^*$ of order $d$. 

A simple calculation shows that $\sigma_i(\Kl_n(q,a)) = \Kl_n(q, i^{n+1} a)$ for all $i\in\F_p^*$. Put $d:= \gcd(n+1,p-1)$ and let $i_d$ be a generator for the unique subgroup of order $d$ of $\F_p^*$. Then $\sigma_{i_d}(\Kl_n(q,a))=\Kl_n(q,i_d^{n+1} a)=\Kl_n(q,a)$ so that $\Kl_n(q,a)\in E_{d}$ for all $a\in\F_q$. On the converse direction we have the following result by Wan \cite{wan:min_pol_and_distinctness}. 
\begin{theorem}\label{wan's_result}
Let $a\in\F_q^*$ be such that $\Tr(a)\not= 0$. Then $\Kl_n(q,a)$ generates the unique intemediate field between $\Q$ and $\Q(\zeta_p)$ of degree $\frac{p-1}{d}$ over $\Q$. 
\end{theorem}

In the present work our aim is to further characterize the field generated by $\Kl_n(q,a)$ also in the case $\Tr(a)=0$. This is possible by using a method we already briefly outlined in \cite{kononen:on_int_values_it}. We consider when $\sigma_i(\Kl_n(q,a))=\Kl_n(q,a)$, ie. when 
\begin{equation}\label{basic_invariance}
  \Kl_n(q,i^{n+1} a) = \Kl_n(q,a). 
\end{equation}

Since the Kloosterman sum is invariant under Frobenius conjugation, ie. $\Kl_n(q,a^{p^s})=\Kl_n(q,a)$, \eqref{basic_invariance} is true if there exists some $s\in\{0,1,\ldots,r-1\}$ such that 
\begin{equation}\label{basic_criteria}
  i^{n+1} a = a^{p^s}. 
\end{equation}
Furthermore there exist some fields for which the values of Kloosterman sums are distinct except Frobenius conjugation ie. that satisfy 
\begin{equation}\label{dist_prop}
 \Kl_n(q,a) = \Kl_n(q,b) \implies b=a^{p^s} \text{ for some } s\in\{ 0,1, \ldots, r-1 \}.
\end{equation}
For these fields \eqref{basic_invariance} holds if and only if \eqref{basic_criteria} is true for some $s$ and we obtain the following theorem which completely characterizes the field generated by $\Kl_n(q,a)$. 
\begin{theorem}\label{main_theorem}
Assume that the field $\F_q$ satisfies \eqref{dist_prop} and denote $R:= \gcd(\tfrac{p-1}{d},r)$. Also, let $\gamma$ be a primitive element for $\F_q$. Then each Kloosterman sum $\Kl_n(q,a)$ for $a\in\F_q^*$ generates some intermediate field between $E_{dR}$ and $E_d$. For $e|R$, $e>1$ we have $\Kl_n(q,a)\in E_{de}$ if and only if there exists a positive divisor $t$ of $\tfrac{r}{e}$ and $u\in\{1,\ldots,e-1\}$ satisfying $\gcd(e,u)=1$ such that $a\in \sqrt[e]{\gamma_t}^u \F_{p^t}^*$ where 
\begin{equation*}
\sqrt[e]{\gamma_t}:= \gamma^{\frac{p^r-1}{e(p^t-1)}}
\end{equation*} 
is an $e$th root of the primitive element $\gamma_t$ of $\F_{p^t}$.  
\end{theorem}

Fischer \cite{fisher:distinctness_of_kloost} has shown that \eqref{dist_prop} holds for all fields satisfying the bound $p>(2(n+1)^{2r}+1)^2$. The referee of the Fischer's paper conjectured that the much stricter bound $p\geq r(n+1)$ should actually hold. Further remarkable progress on this question has been made by Wan in the previously mentioned work \cite{wan:min_pol_and_distinctness}. Among other things Wan showed that \eqref{dist_prop} holds with the bound $p\geq(r-1)(n+1)+2$ provided the prime $p$ does not divide any of the terms in a certain rather complicated finite sequence. Especially it follows from Wan's results that referee's conjecture holds in the case of classical Kloosterman sums ($n=1$) for many small values of $r$ except possibly some small number of exceptional values of $p$. In particular, it holds for all $p$ in the case $n=1$, $r\leq 4$. 

Even if the field $\F_q$ does not satisfy \eqref{dist_prop} we still have the next partial result. 
\begin{theorem}\label{main_result2}
Let $e|\gcd(\tfrac{p-1}{d},r)$, $e>1$. We have $\Kl_n(q,a)\in E_{de}$ if there exists a positive divisor $t$ of $\tfrac{r}{e}$ and $u\in\{1,\ldots,e-1\}$ satisfying $\gcd(e,u)=1$ such that $a\in \sqrt[e]{\gamma_t}^u \F_{p^t}^*$ (using the same notation as in Theorem \ref{main_theorem}).
\end{theorem}

\section{The proof of Theorems \ref{main_theorem} and \ref{main_result2}}

We need the following elementary fact.

\begin{lemma}\label{residue_lemma}
Let $m|n$. Then $x+n\Z\mapsto x+m\Z$ defines a surjective homomorphism $(\Z/n\Z)^* \to (\Z/m\Z)^*$. 
\begin{proof}
Write $n=MQ$ where $\gcd(Q,m)=1$ and $M$ is the product of all prime factors $l$ of $n$ such that $l|m$ (counting multiplicities). Let $y+m\Z\in (\Z/m\Z)^*$. We shall show that at least one of $y, y+m, \ldots, y+(Q-1)m$ is relatively prime to $n$. 

Put $k_i:= \gcd(n, y+im)$ for every $i\in\Z$. We have $\gcd(m, y+im)=\gcd(m,y)=1$ and therefore also $\gcd(M,y+im)=1$. It follows that $k_i=\gcd(Q,y+im)$ for every $i$. Notice that the elements $(y+im)+Q\Z$ for $i=0,1, \ldots, Q-1$ form a coset of the additive subgroup of $\Z/Q\Z$ generated by $m$. Since $\gcd(m,Q)=1$ this subgroup is actually whole $\Z/Q\Z$. Consequently, there exists an $i\in\{0,1,\ldots,Q-1\}$ such that $k_i=\gcd(Q,y+im)=1$. 
\end{proof}
\end{lemma}

Denote $\gamma_t=\gamma^{(p^r-1)/(p^t-1)}$ whenever $t|r$ so that $\gamma_t$ is a primitive element for $\F_{p^t}$. Assume that $e|\tfrac{p-1}{d}$ and $a\in\F_q^*$. Our idea is to determine when  $\Kl_n(q,a)\in E_{de}$. Obviously this is the case if and only if \eqref{basic_invariance} holds   
for a generator $i$ of the order $ed$ subgroup of $\F_p^*$. We shall use $i=\gamma_1^{(p-1)/ed}$. Obviously \eqref{basic_invariance} is satisfied if there exists an $s\in\{0,1,\ldots,r-1\}$ such that \eqref{basic_criteria} holds. Moreover, if we assume \eqref{dist_prop} (that is in the case of Theorem \ref{main_theorem}) we find that \eqref{basic_invariance} is satisfied if and only if there exists an $s\in\{0,1,\ldots,r-1\}$ such that \eqref{basic_criteria} holds.

Let us now try to find out all $a$ and $s$ such that \eqref{basic_criteria} is satisfied for a fixed $e$. For $s=0$ and any $a\in\F_q^*$ we see that \eqref{basic_criteria} is true if and only if $i^{n+1}=1$. This corresponds to the trivial case $e=1$. From now on we shall assume $e>1$ and $s>0$. Put $a=\gamma^x$. Then \eqref{basic_criteria} is equivalent to 
\begin{equation}\label{congruence1}
 (p^s-1) x \equiv \frac{n+1}{d} \frac{p^r-1}{e}  \pmod{p^r-1}. 
\end{equation} 
Note that $\gcd(\tfrac{n+1}{d},\tfrac{p-1}{d})=1$ by definition of $d$ and therefore also $\gcd(e,\tfrac{n+1}{d})=1$.

We have $\gcd(p^s-1,p^r-1)=p^t-1$ where $t:= t_s:= \gcd(r,s)$. Put $r=r't$ and $s=s't$. Now \eqref{congruence1} is solvable if and only if $(p^t-1)|\tfrac{n+1}{d}\tfrac{p^r-1}{e}$ ie. if and only if $e|\tfrac{n+1}{d}\tfrac{p^r-1}{p^t-1}$. Because $\gcd(e,\tfrac{n+1}{d})=1$ this is equivalent to $e|\gcd(\tfrac{p-1}{d}, \tfrac{p^r-1}{p^t-1})$. Notice that \begin{equation*}
\frac{p^r-1}{p^t-1}=(p^t-1)(p^{t(r'-2)}+2p^{t(r'-3)}+\cdots+(r'-2)p^t+(r'-1)) + r'. 
\end{equation*} 
Since $\tfrac{p-1}{d}|(p^t-1)$ we must have $\gcd(\tfrac{p-1}{d},\tfrac{p^r-1}{p^t-1})=\gcd(\tfrac{p-1}{d},r')$. Thus, for a fixed $e$ and $s$ satisfying $\gcd(r,s)=t$, \eqref{congruence1} is solvable if and only if $e|\gcd(\tfrac{p-1}{d},\tfrac{r}{t})$. 

Let us now assume that \eqref{congruence1} is solvable for a given $s$ and let us find all solutions $a$ for this $s$. The congruence \eqref{congruence1} will then be equivalent to 
\begin{equation}\label{congruence2}
 \frac{p^s-1}{p^t-1}x \equiv \frac{n+1}{d}\frac{p^r-1}{e(p^t-1)}\pmod{\frac{p^r-1}{p^t-1}}.
\end{equation}
Denote $C_{s,t}:= \tfrac{p^s-1}{p^t-1}$ for all $s,t$. Assume that $C_{s,t} C_{s,t}^* \equiv 1 \pmod{C_{r,t}}$. Then we can write the solutions for \eqref{congruence2} in the form $x\equiv \tfrac{n+1}{d} C_{s,t}^* \tfrac{C_{r,t}}{e} \pmod{C_{r,t}}$. We have $p\equiv 1 \pmod e$ and therefore also 
\begin{equation*}
C_{s,t} = \tfrac{p^{ts'}-1}{p^t-1} = p^{t(s'-1)}+\cdots + p^{t}+1 \equiv s' \pmod{e}.
\end{equation*} 
Using this information the solution for \eqref{congruence2} comes 
\begin{equation*} 
x\equiv \tfrac{n+1}{d} \tfrac{s^* C_{r,t}}{e} \pmod{C_{r,t}}
\end{equation*} 
where $s^*$ is such that $s' s^* \equiv 1 \pmod e$.

Notice that if \eqref{congruence1} is solvable for some $s$ then it will also be solvable (with the same $e$) for all other $s$ satisfying $\gcd(r,s)=t$. There are $\Euler(\tfrac{r}{t})$ such values of $s$, namely $s=s't$, where $1\leq s' <\tfrac{r}{t}$ and $\gcd(s',r')=1$. Since $e|\tfrac{r}{t}$ it follows from Lemma \ref{residue_lemma} that $s'\mod e$ and consequently also $s^*\tfrac{n+1}{e}$ will run through the reduced residue system modulo $e$ (possibly several times) when $s$ goes through all the above values. Thus all the solutions $x$ for all $s$ giving rise to fixed $t|r$ and for this fixed $e$ will be 
\begin{equation*}
 x = u\frac{p^r-1}{e(p^t-1)} + v\frac{p^r-1}{p^t-1}
\end{equation*}
where $u\in\{1,2,\ldots, e-1\}, \gcd(e,u)=1$ and $v\in\{0,1,\ldots, p^t-2\}$. Obviously we can write these solutions in the claimed form $a=\gamma^x\in \sqrt[e]{\gamma_t}^u \F_{p^t}^*$. 

From the preceeding discussion it follows that there exists $t$ such that \eqref{congruence1} is solvable for the corresponding $s$ if and only if $e|R$. In particular, under the assumption \eqref{dist_prop} this means that $\Kl_n(q,a)$ generates some intermediate field of $E_{dR}$ and $E_d$ for all $a\in\F_q^*$. Furthermore, if solutions exist then the set of all $t$ for which \eqref{congruence1} is solvable consists positive divisors $t|\tfrac{r}{e}$ (assuming $e>1$). Gathering all the information together we have the claimed results.

\section{Further observations}

We note that in Theorems \ref{main_theorem} and \ref{main_result2} the condition denoted by (i) below can be written in the equivalent form (ii).
\begin{itemize}
\item[(i)  ] $a \in \sqrt[e]{\gamma_t}^u \F_{p^t}^*$, for some $u\in\{0,1,\ldots,e-1\}$ relatively prime to $e$. 
\item[(ii) ] $k=e$ is the smallest positive exponent such that $a^k\in\F_{p^t}^*$. 
\end{itemize}

It follows from the Wan's result (Theorem \ref{wan's_result}) that $\Tr(a)=0$ whenever $\Kl_n(q,a)\in E_{de}$ for $e>1$. We can, however, also see this directly for the values of $a$ obtained from the preceeding theorems. Assume $e>1$, $e|\gcd(\tfrac{p-1}{d},r)$ and $t|\tfrac{r}{e}$. Let $a\in\F_q^*$ be such that $k=e$ is the smallest exponent for which $a^k\in\F_{p^t}$. We have already seen in the previous section that $p\equiv 1 \pmod e$ implies $\tfrac{p^{et}-1}{p^t-1}\equiv e \pmod e$. It follows that $\sqrt[e]{\gamma_t}\in\F_{p^{et}}$ and thus also $a\in\F_{p^{et}}$. We have 
\begin{equation*}
  \Tr_{\F_{p^{et}}/\F_{p^t}}(a) = a \left(1+a^{p^t-1}+\cdots + a^{p^{t(e-1)}-1}\right).
\end{equation*}
Since $a^e\in\F_{p^t}$ and $e|(p-1)$ we see that $(1+a^{p^t-1}+\cdots + a^{p^{t(e-1)}-1})\in\F_{p^t}$. By the properties of the trace function also $\Tr_{\F_{p^{et}}/\F_{p^t}}(a)\in\F_{p^t}$ and by assumptions $a\not\in\F_{p^t}$. This is possible if and only if 
\begin{equation*}
\Tr_{\F_{p^{et}}/\F_{p^t}}(a) = 0.
\end{equation*} 
From this sharper result it follows, by using the transitivity of the trace, that $\Tr(a)=0$. 

As a special case of our results we can characterize some -- all in the case of \eqref{dist_prop} -- points where the Kloosterman sum obtains a rational value.
\begin{corollary}
Assume that $\tfrac{p-1}{d}|r$. We have $\Kl_n(q,a)\in \Q$ if there exists a positive divisor $t$ of $\tfrac{dr}{p-1}$ such that $k=\tfrac{p-1}{d}$ is the smallest positive exponent for which $a^k\in\F_{p^t}$.

Furthermore, if $\F_q$ satisfies \eqref{dist_prop} there are no other values $a\in\F_q^*$ besides the above for which $\Kl_n(q,a)\in\Q$. In particular, in the case $\tfrac{p-1}{d}\nmid r$ there exists no $a\in\F_q^*$ for which $\Kl_n(q,a)\in \Q$. 
\end{corollary}


\end{document}